\documentclass[review]{elsarticle}

\usepackage{lineno,hyperref}
\modulolinenumbers[5]
\usepackage{amsfonts}
\usepackage{amsmath}
\usepackage{amsthm}
\newcommand{\lap}{\mbox{$\bigtriangleup$}}
\newtheorem{thm}{Theorem}[section]

\theoremstyle{definition}

\theoremstyle{remark}
\newtheorem{rem}{Remark}[section]

\newcommand{\R}{\mathbb{R}}

\journal{Journal of \LaTeX\ Templates}

\begin{document}
\title{\bf Symmetric properties for Choquard equations involving fully nonlinear  nonlocal operator }
\author{Pengyan Wang$^a$$^b$, Li Chen$^b$, Pengcheng Niu$^a$ \thanks{The corresponding author's email address: chen@math.uni-mannheim.de(L. Chen)}\\
\small{$^a$ Department of Applied Mathematics, Northwestern Polytechnical
University,}\\ \small{ Xi'an, Shaanxi, 710129, P. R. China}\\
\small{$^b$ School of business informatics and mathematics, }\\
\small{Universitat of Mannheim, 68131, Mannheim, Germany }\\
\small e-mail: wangpy119@126.com, chen@math.uni-mannheim.de,pengchengniu@nwpu.edu.cn}
\date{}
\maketitle{\bf Abstract}.
In this paper, the positive solutions to Choquard equation involving fully nonlinear nonlocal operator are shown to be symmetric and monotone by using the  moving plane method which has been introduced by Chen, Li and Li in 2015. The key ingredients are to obtain the ``narrow region principle'' and ``decay at infinity'' for the corresponding problems. Similar ideas can be easily applied to various nonlocal problems with more general nonlinearities.

\noindent{{\bf Keywords}:  Fully
nonlinear nonlocal Choquard equation, method of moving planes, narrow region princple,  decay at infinity.

\section{Introduction}
In this paper, we study the following Choquard equation    involving fully nonlinear nonlocal operator:
\begin{equation}
\label{eq:m1} \left\{\begin{array}{ll}
{\cal F} _{\alpha}(u(x))+\omega u(x) =C_{n,2s} (|x|^{2s-n}\ast u^q(x)) u^r(x),  ~& x\in \R^n, \\
u (x)>0,~ &  x\in \R^n,
\end{array}
\right.
\end{equation}
where the operator $\cal F_\alpha$ with $0<\alpha<2$ is given by
$$
{\cal{F}}_{\alpha}(u(x))
= C_{n,\alpha} \mbox{ P.V. } \int_{ {\R}^n} \frac{F (u (x)-u (y))}{|x-y|^{n+\alpha }} dy=C_{n,\alpha} \lim_{\epsilon \rightarrow 0} \int_{ {\R}^n\setminus B_{\epsilon}(x)} \frac{F(u(x)-u(y))}{|x-y|^{n+\alpha}} dy.
$$
In this formulation $\mbox{ P.V. }$ stands for the Cauchy principal value of the integral, $F$ is a given local Lipschitz continuous function defined on $\R$, with  $F(0)=0$ and $F'\geq c_0>0$, and $C_{n,\alpha}=(2-\alpha)\alpha2^{\alpha-2}\frac{\Gamma(\frac{n+\alpha}{2})}{ {\pi}^{\frac{n}{2}}\Gamma(\frac{4-\alpha}{2})}$.

The Choquard equation \eqref{eq:m1} is considered in the case that $\omega\geq 0$, $0<s<1$ and $1<r,q<\infty$. $C_{n,2s}$ has the same representation to $C_{n,\alpha}$.

The fully nonlinear nonlocal operator has been introduced by
Caffarelli and Silvestre in \cite{CS1}.
Define
\begin{eqnarray*}
{\cal L}_{\alpha  }(\R^n)=\Big\{ u : F(u)\in L^1_{loc}(\R^n), \,   \int_{ {R}^{n}} \frac{|F (1+u (x))|}{1+|x|^{n+\alpha }}dx <\infty \Big\},
\end{eqnarray*}
then it is easy to see that for $u \in C^{1,1}_{loc} \cap {\cal L}_{\alpha },$ ${\cal F} _{\alpha }(u)$ is well defined.
In the case that $F (\cdot)$ is linear function, ${\cal F}_{\alpha}$ becomes  the  usual  fractional Laplacian $(-\Delta)^{\frac{\alpha}{2}}$, $0<\alpha<2$. We refer to \cite{CLM} for the definition and further properties of fractional Laplacian.   When $F (t)=|t|^{p-2}t,$ $\alpha =\theta p$, $0<\theta<1$  and $1<p<\infty,~$   ${\cal F}_{\alpha}$ becomes  fractional $p$-Laplacian $(-\Delta)_p^{\theta}$.
The nonlocality of the fractional Laplacian makes it difficult to investigate. To our knowledge, the extension
method   \cite{CS} and the integral equation approach \cite{CLO} have been employed successfully to study equations involving the fractional Laplacian. However, both methods can not be directly used to handle equations involving fully nonlinear nonlocal operator. In 2015, Chen, Li and Li \cite{CLL}   developed a new technique (the direct method of moving planes) that can be applied for problems with fractional Laplace operator. It is very effective in dealing with equations involving fully nonlinear nonlocal operators or uniformly elliptic nonlocal operators, for example the results in \cite{CWD,CLLg,WN,WW,WY}  and the references therein.

Back to equation \eqref{eq:m1}, when $F(t)=t$, $s=1$, $r=2$, $q=1$, $\omega>0$ and $\alpha= 2$, it is reduced to the well known Choquard or   nonlinear Hartree equation
\begin{equation}\label{eq:mmw1}
-\lap u(x) +\omega u(x) = (|x|^{2 -n}\ast u^2(x)) u(x) ,  ~  x\in \R^n.
 \end{equation}
Equation \eqref{eq:mmw1} in three dimension has strong background in quantum mechanics.
It is the one particle mean field approximation for many particle interacting Coulomb system, which has be proved rigorously in the last decades by many authors. We refer only a few of them, \cite{LNR,LSSY,LY}.   If $u$ solves \eqref{eq:mmw1}, then the function $\psi$ defined
by $\psi(x,t)=e^{i\omega t}u(x)(\omega>0)$ is a solitary wave solution of the focusing time-dependent Hartree equation
\begin{equation}\label{eq:mw3}
i\psi_t+ \lap \psi(x) +    (|x|^{2 -n}\ast \psi^2) \psi = 0,  ~  x\in \R^n, ~t>0.
  \end{equation}
The rigorous derivation of equation \eqref{eq:mw3} via mean field limit of many body Schr\"odinger equation with two body interaction has been extensively studied by several research groups in mathematical physics. For example the convergence results in \cite{EY} by using the BBGKY hierarchy, the convergence rate estimate in \cite{CLi, CLiS,P,RS} to name a few. And the uniqueness and qualitative results of the elliptic system   in  \cite{BTWW} have aroused our interest. We also refer readers to related Choquard and Schr\"{o}dinger equations which involve local and nonlocal operators, \cite{RK,RK1,MpZ} for the existence of   solutions, \cite{MZl} for the classification of positive solutions, \cite{MS} for the qualitative properties and decay asymptotics, \cite{MV} for a review, and references therein.

In this article, we study the symmetry of positive solutions for the Choquard type equation \eqref{eq:m1}.  The symmetry of solutions plays important roles in the  qualitative analysis of the solutions for partial differential equations.
The method of moving planes introduced by Alexanderov \cite{A} in the early 1950s is a powerful tool in obtaining the symmetry of solutions to
elliptic equations and systems. It has been further developed by Serrin \cite{S}, Gidas, Ni and Nirenberg \cite{GNN}, Caffarelli, Gidas and Spruck \cite{CGS},  Chen and Li \cite{CL},   Li and Zhu \cite{LZh}, Li \cite{L4},
Lin \cite{Lin}, Chen, Li and Ou \cite{CLO}, Chen, Li and Li \cite{CLL} and many others. Here, we   apply the direct moving plane method which  introduced by  Chen, Li and Li \cite{CLL} to equations involving two different types of nonlocal operators.

Throughout the paper, let \begin{equation}\label{eq:m2}\aligned
v(x)=C_{n,2s}(|x|^{2s-n}\ast u^q(x))=C_{n,2s}\int_{\R^n}\frac{u^q(y)}{|x-y|^{n-2s}}dy.
\endaligned\end{equation}
Note the Green's function of $(-\lap)^s$ in $\R^{n}$ is $ \frac{C_{n,2s}}{|x-y|^{n-2s}}$, for $v\in C^{1,1}_{loc}\cap{L}_{ 2s}(\R^n)$, where $L_{2s}=\{u \in L^1_{loc}(\mathbb R^n):\R^n\rightarrow \mathbb R\mid \int_{\R^n}\frac{|u(x)|}{1+|x|^{n+2s}}dx<\infty\}$, the following equivalent formulation holds.
$$
(-\lap)^sv(x)=u^q(x),x\in \R^n.
$$
Hence, \eqref{eq:m1} is equivalent to
\begin{equation}\label{eq:m3}
\left\{\begin{array}{ll}
  {\cal F}_{\alpha}(u(x))+\omega u(x)=v(x)u^r(x),&\quad  x\in \R^n,\\
  (-\lap)^sv(x)=u^q(x),  &\quad x\in \R^n,\\
  u(x)>0,v(x)>0, &\quad x\in \R^n.\end{array}
\right.
\end{equation}
Therefore, the investigation of \eqref{eq:m1} is reduced to study \eqref{eq:m3}.

The main result of this paper is
\begin{thm}\label{thmm}
	Suppose that $u \in  C^{1,1}_{loc}(\R^n) \cap {\cal L}_\alpha(\R^n)$ be   positive solution  of equation \eqref{eq:m1} with $1<r,q <\infty$  and satisfy
	\begin{equation}\label{eq:mudecay}
	u(x)\sim \frac{1}{|x|^\gamma}  \quad   \mbox{for} ~|x| ~\mbox{sufficiently large},
	\end{equation}
	where $\gamma$ satisfy
	\begin{equation}\label{eq:gamma}\aligned
	\frac{n}{q}>\gamma>\max\{\frac{\alpha}{r},\frac{2s}{q-1}, \frac{2s+\alpha}{r-1+q}\}.
	\endaligned\end{equation}
	Then $u$  is radially symmetric and  monotone decreasing  about some point in  $\R^n$.
\end{thm}
\begin{rem}
	The assumption \eqref{eq:gamma} for the existence of $\gamma$ implies  $n>q\max\{\frac{\alpha}{r},\\
\frac{2s}{q-1}, \frac{2s+\alpha}{r-1+q}\} $ is required.
\end{rem}
Due to the equivalence of problems \eqref{eq:m1} and \eqref{eq:m3}, we only need to prove the following theorem for \eqref{eq:m3}.
\begin{thm}\label{thmm1}
Suppose that $u \in  C^{1,1}_{loc}(\R^n) \cap {\cal L}_\alpha(\R^n)$ and $v \in  C^{1,1}_{loc}(\R^n) \cap {L}_{2s}(\R^n)$ be   positive solutions  of equation \eqref{eq:m3} with $1<r,q <\infty$  and satisfy
 \begin{equation}\label{eq:m4}
u(x)\sim \frac{1}{|x|^\gamma} ,\quad  v(x)\sim \frac{1}{|x|^\tau} \quad   \mbox{for} ~|x| ~\mbox{sufficiently large},
\end{equation}
where $\gamma>0,\tau>0$ satisfy
\begin{equation}\label{eq:m5}\aligned
 \alpha<\min \{\gamma(r-1)+\tau,r\gamma\}~\mbox{and}~ 2s<  \gamma(q-1).
\endaligned\end{equation}
Then $u$ and $v$ are radially symmetric and  monotone decreasing  about some point in  $\R^n$.
 \end{thm}

\begin{rem}
	From the assumptions \eqref{eq:mudecay} and \eqref{eq:gamma} in theorem \ref{thmm}, one can obtain that the assumptions \eqref{eq:m4} and \eqref{eq:m5} in theorem \ref{thmm1} are valid with $\tau=q\gamma-2s<n-2s$ by using the Lemma 2.1 in \cite{CLZ}.
\end{rem}

\begin{rem}  Since the Kelvin transform is not applicable for the fully nonlinear nonlocal equation, further
assumptions on the behaviors of $u$ and $v$ at infinity are needed.  In case that $\omega=0$, $F (t)=|t|^{p-2}t$, and $\alpha =\theta p$, Theorem \ref{thmm1} recovers the result of Choquard equation involving fractional $p$-Laplacian in \cite{MZ}.
  \end{rem}

The paper is organized as follows. In Section 2, we  establish the ``narrow region principle''  and ``decay at infinity''.   Section 3 is devoted to the  proof  of   \ref{thmm1}   by using the method of moving planes. As already has been remarked, the results of Theorems \ref{thmm} can be obtained directly.

Throughout the paper,    $C$ will be positive constants  which can  be different from line to line and only the relevant
dependence is specified.

\section{``Narrow region principle''  and ``decay at infinity'' }

Through out this section, let  $u \in  C^{1,1}_{loc}(\R^n) \cap {\cal L}_\alpha(\R^n)$ and $v \in  C^{1,1}_{loc}(\R^n) \cap {L}_{2s}(\R^n)$ be positive solutions  of equation \eqref{eq:m3} which satisfy the assumptions of theorem \ref{thmm1}.

For simplicity,  we list some notations that will used frequently:  for $\lambda \in \R, $ denote $x=(x_1,x'),$
$x^\lambda=(2\lambda-x_1,x')$  and
$$
T_\lambda=\{x\in \R^n| x_1=\lambda \}, \quad \Sigma_\lambda=\{x \in \R^n\mid x_1<\lambda\}, \quad \tilde \Sigma_\lambda =\{x\in \R^n\mid x_1>\lambda\}.
$$
We denote
\begin{eqnarray*}
u _\lambda(x )=u (x^\lambda ), && v _\lambda(x )=v (x^\lambda ),\\
U_\lambda (x )=u _\lambda(x )
-u ( x ), &&V_\lambda (x )=v _\lambda(x  )
-v ( x ).
\end{eqnarray*}
We call that a function $U_\lambda(x )$ is  a  $\lambda$ antisymmetric function   if and only if
\begin{equation}\label{eq:j2} U_\lambda(x_1,x_2,\cdots,x_n )=-U_\lambda(2\lambda-x_1,x_2,\cdots,x_n ). \end{equation}
Obviously, $U_\lambda(x)$ and $V_\lambda(x)$  are antisymmetric functions.

\begin{thm}\label{thmm2}(narrow region principle  )
Let $\Omega$ be a region contained in
 $$\{x |  \lambda-l<x_1<\lambda  \}\subset \Sigma_\lambda$$
 with small $l$.  Suppose  that $U_\lambda(x)\in C_{loc}^{1,1}(\R^n)\cap{\cal L}_\alpha$ and  $ V_\lambda(x)\in C_{loc}^{1,1}(\R^n)\cap{ L}_{2s} $   are lower semi-continuous on $\bar{\Omega}$,
  and satisfy
\begin{equation}\label{eq:m6}
\left\{\begin{array}{ll}
{\cal F}_\alpha (u_\lambda(x))-{\cal F}_\alpha (u(x))+c_{1}(x)U_\lambda(x)+c_{2}(x)V_\lambda(x)\geq0, & x \in \Omega, \\[2mm]
(-\lap)^s V_\lambda(x)  + c_{3}(x)U_\lambda(x) \geq0, & x \in \Omega,\\[2mm]
U_\lambda(x), V_\lambda(x)\geq 0, & x \in \Sigma_\lambda\backslash\Omega,\\[2mm]
U_\lambda(x^\lambda)=-U_\lambda(x),
V_\lambda(x^\lambda)=-V_\lambda(x), & x \in \Sigma_\lambda.
\end{array}
\right.
\end{equation}
where  $c_{i}(x)(i=1,2,3)$ are bounded from below in $\Omega$,~ $c_{2}(x)<0,c_{3}(x)<0$, then the following statements hold.
\begin{enumerate}
\item[(i)]   If $\Omega$ is bounded, then  for sufficiently small $l$,
\begin{equation}\label{eq:mb4}
U_\lambda(x),V_\lambda(x) \geq0 , \forall~ x\in \Omega;
\end{equation}
\item[(ii)]  if $\Omega$ is unbounded, then the conclusion \eqref{eq:mb4} still holds under the conditions
\begin{equation}\label{eq:mbbb4}
\underset{|x|\rightarrow \infty}{\underline{\lim}}U_\lambda(x),~\underset{|x|\rightarrow \infty}{\underline{\lim}}V_\lambda(x)\geq0;
\end{equation}
\item[(iii)]{(Strong maximum priciple)} If  \eqref{eq:mb4} holds  and $U_\lambda(x)$ or $V_\lambda(x)$ attains 0 somewhere in $\Omega$, then
\begin{equation}
\label{eq:mb6}U_\lambda(x)=V_\lambda(x)\equiv 0,~ x\in \R^n.
\end{equation}
\end{enumerate}
\end{thm}

\begin{rem}  If $U_\lambda(x)$ or $V_\lambda(x)$ is positive at some point in $\Omega$, then
it follows by (iii)  that
$$U_\lambda(x)>0,~V_\lambda(x)> 0,~\forall~ x\in \Omega.$$
As we can see from the proof of (iii) later, $\Omega$ can be bounded or unbounded and does not need to be narrow. \end{rem}

\begin{proof} (i)
Suppose that \eqref{eq:mb4} does not hold, without loss of generality, we assume $V_\lambda(x)$ is negative at some point in $\Omega$; then the lower semi-continuity of $V_\lambda(x)$ on $\bar \Omega$ implies that  there exists  $\bar x$ such that
$$V_\lambda(\bar x)=\underset{\Omega}{\min} V_\lambda(x)<0,$$
and $\bar{x}$ is in the interior of $\Omega$ from the  condition \eqref{eq:m6}.
Similar to  the calculation in \cite{CLL}, we can derive that
\begin{equation}\label{eq:m99}(-\lap)^sV_\lambda(\bar x)  \leq C_{n,2s} \int_{\Sigma_\lambda}\frac{2V_\lambda(\bar x)}{|\bar x-y|^{n+ 2s}}dy \leq \frac{C V_\lambda(\bar  x)}{2l^{ 2s}}<0.\end{equation}
From the second inequality of \eqref{eq:m6} one has
\begin{equation}\label{eq:m20}0\leq(-\lap)^sV_\lambda(\bar x)  +c_{3}(\bar x)U_\lambda(\bar  x)\leq \frac{C V_\lambda(\bar  x)}{2l^{2s}}+c_{3}(\bar x)U_\lambda(\bar  x) .\end{equation}
Since $c_3(\bar x)<0,$ we can drive that
\begin{equation}\label{eq:mb8}
U_\lambda(\bar x)< 0~\mbox{and}~V_\lambda(\bar x)\geq -C c_{3}(\bar x)l^{ 2s} U_\lambda(\bar x),
\end{equation}
which implies that there exists $\tilde{x}\in \Omega $ such that
$$U_\lambda(\tilde{x})=\underset{\Omega} {\min} U_\lambda(x) <0.$$
By   the  expression of ${\cal F}_\alpha$ and  (2.2) in \cite{WY}, we have
 \begin{eqnarray}\label{eq:mbb20}
{\cal F}_\alpha (u_\lambda(\tilde x))-{\cal F}_\alpha (u(\tilde x))
  \leq  2C_{n,\alpha} C U_\lambda(\tilde x)\int_{\Sigma_\lambda}\frac{1}{|\tilde x -y^\lambda|^{n+\alpha}}dy.
\end{eqnarray}
  Then for $0<r<\min\{l-\tilde x_1\}$,
\begin{equation}\label{eq:mb666}\aligned
\int_{\Sigma_\lambda}\frac{1}{|\tilde x -y^\lambda|^{n+\alpha}}dy
\geq & \int_{ B_r(\tilde x)}\frac{1}{|\tilde x -y^\lambda|^{n+\alpha}}dy\\
\geq &\int_{ B_r(\tilde x^\lambda)}\frac{1}{|\tilde x -y|^{n+\alpha}}dy\geq\frac{w_n}{3r^\alpha}\geq\frac{C}{l^\alpha}.\endaligned
\end{equation}
Plug it into \eqref{eq:mbb20}, it follows that
\begin{equation}\label{eq:m10}{\cal F}_\alpha (u_\lambda(\tilde x))-{\cal F}_\alpha (u(\tilde x))\leq \frac{C U_\lambda(\tilde x)}{l^\alpha}<0.\end{equation}
which together with \eqref{eq:m6}, for $l$ sufficiently small, implies that
\begin{equation}\label{eq:mb7}{\cal F}_\alpha (u_\lambda(\tilde x))-{\cal F}_\alpha (u(\tilde x))+c_{ 1}(\tilde x)U_\lambda(\tilde x)\leq (\frac{C}{l^\alpha}+c_{ 1}(\tilde x))U_\lambda(\tilde x)\leq \frac{C}{ l^\alpha}U_\lambda(\tilde x)<0.\end{equation}
Combining \eqref{eq:m6}, \eqref{eq:mb8} and \eqref{eq:mb7}, for $l$ sufficiently small, we have
$$\aligned
0&\leq {\cal F}_\alpha (u_\lambda(\tilde x))-{\cal F}_\alpha (u(\tilde x))+c_{ 1}(\tilde x)U_\lambda(\tilde x)+c_{ 2}(\tilde x)V_\lambda(\tilde x)\\
&\leq\frac{CU_\lambda(\tilde x)}{ l^\alpha}+c_{2 }(\tilde x)V_\lambda(\bar x)\\
&\leq\frac{CU_\lambda(\tilde x)}{l^\alpha}-Cc_{2 }(\tilde x)c_{3}(\bar x)l^{2s} U_\lambda(\bar x)\\
&\leq\frac{CU_\lambda(\tilde x)}{l^\alpha}-Cc_{2 }(\tilde x)c_{3}(\bar x)l^{2s} U_\lambda(\tilde x)\\
&\leq C \frac{U_\lambda(\tilde x)}{l^\alpha}(1-c_{ 2}(\tilde x)c_{3}(\bar x)l^{\alpha+2s})<0.
\endaligned$$
This contradiction shows that \eqref{eq:mb4} must be true.

(ii)  If $\Omega$ is unbounded, then  \eqref{eq:mbbb4} guarantees that the negative minimum of $U_\lambda$ and $V_\lambda$ must be attained at some point $\tilde x$ and $\bar x $, respectively. Then one can follow the same discussion as the case of (i) to arrive at a contradiction.

(iii) To prove  \eqref{eq:mb6}, without loss of generality,  we suppose that there exists $z\in \Omega$ such that
$$V_\lambda(z)=0.$$
Then  $\frac{1}{|x-y|}>\frac{1}{|x-y^\lambda|},~\forall x,y\in \Sigma_\lambda$ and
\begin{equation}\label{eq:mbb10}\aligned
 &(-\lap)^s V_\lambda(z) =C_{n,2s} \mbox{ P.V. }\int_{\R^{n}}\frac{ -V_\lambda(y) }{|z -y|^{n+2s}}dy\\
& \leq  C_{n,2s}\mbox{ P.V. }  \int_{\Sigma_\lambda}V_\lambda(y)(\frac{1}{|z -y^\lambda|^{n+2s}}-\frac{1}{|z -y|^{n+2s}})dy.
 \endaligned
\end{equation}
If $V_\lambda(y)$ is not identically equals to zero in $\Sigma_\lambda$, then \eqref{eq:mbb10} implies that
\begin{equation}
\label{eq:m11}(-\lap)^sV_\lambda(z)<0.
\end{equation}
Combining \eqref{eq:m11} with the  second inequality of \eqref{eq:m6}, we get
$$ U_\lambda(z)<0.$$
This is a contradiction with \eqref{eq:mb4}. Hence $V_\lambda(x)$ must be identically $0$ in $\Sigma_\lambda$.
Since $$V_\lambda(x^\lambda)=-V_\lambda(x),~\forall~x\in \Sigma_\lambda,$$ it yields that
\begin{equation}
\label{eq:m16}V_\lambda(x)\equiv0,~\forall~x\in {R}^n.
\end{equation}
Again from the second equation of \eqref{eq:m6}, \eqref{eq:mbb10} and \eqref{eq:m16}, we know that $$U_\lambda (x)\leq 0,~\forall~x\in \Omega.$$  From \eqref{eq:mb4},
it must hold that
\begin{equation}\label{eq:m17}U_\lambda(x)=0,~\forall~x\in \Omega.\end{equation}
Next we prove $U_\lambda(x)=0,~\forall~x\in \R^n\backslash \Omega.$ If not, we have $U_\lambda(x)\not\equiv  0,~\forall~x\in \R^n\backslash \Omega$.

Now $\forall \tilde x\in \Omega$, it followes from \eqref{eq:m17} that $U_\lambda (\tilde x)=0$.
Therefore one can deduce from  \eqref{eq:m16} and \eqref{eq:m17} that
\begin{eqnarray}
\nonumber &&{\cal F}_\alpha (u_\lambda(\tilde x))-{\cal F}_\alpha (u(\tilde x)) +c_1(x)U_\lambda(\tilde x) +c_2(\tilde x)V_\lambda(\tilde x) \\
\nonumber&=&C_{n,\alpha} \mbox{ P.V. }\int_{\R^{n}}\frac{F(u_\lambda (\tilde x)-u_\lambda (y))-F(u (\tilde x)-u(y))}{|\tilde x -y|^{n+\alpha}}dy\\
\nonumber &=&C_{n,\alpha} \mbox{ P.V. }\int_{\Sigma_\lambda}\frac{F(u_\lambda (\tilde x)-u_\lambda (y))-F(u (\tilde x)-u(y))}{|\tilde x -y|^{n+\alpha}}dy\\
\nonumber &&+C_{n,\alpha} \mbox{ P.V. }\int_{\Sigma_\lambda}\frac{F(u_\lambda (\tilde x)-u (y))-F(u (\tilde x)-u_\lambda(y))}{|\tilde x -y^\lambda|^{n+\alpha}}dy \\
\nonumber &=&C_{n,\alpha} \mbox{ P.V. }\int_{\Sigma_\lambda}\Big[F(u_\lambda (\tilde x)-u_\lambda (y))-F(u (\tilde x)-u(y))\Big]\\
\nonumber &&\hspace{4cm}\cdot \Big(\frac{1}{|\tilde x -y|^{n+\alpha}}-\frac{1}{|\tilde x -y^\lambda|^{n+\alpha}}\Big)dy\\
\nonumber &&+C_{n,\alpha} \mbox{ P.V. }\int_{\Sigma_\lambda}\Big[\frac{F(u_\lambda (\tilde x)-u (y))-F(u (\tilde x)-u_\lambda(y))}{|\tilde x -y^\lambda|^{n+\alpha}}\\
\nonumber &&\hspace{4cm}+\frac{F(u_\lambda (\tilde x)-u_\lambda (y))-F(u (\tilde x)-u(y))}{|\tilde x -y^\lambda|^{n+\alpha}}\Big]dy\\
\nonumber & =&C_{n,\alpha} F'(\cdot)\int_{\Sigma_\lambda} (U_\lambda(\tilde x)-U_\lambda(y))(\frac{1}{|\tilde x -y|^{n+\alpha}}-\frac{1}{|\tilde x -y^\lambda|^{n+\alpha}})dy\\
\nonumber &&+C_{n,\alpha} F'(\cdot)\int_{\Sigma_\lambda}\frac{2U_\lambda(\tilde x)}{|\tilde x -y^\lambda|^{n+\alpha}}dy\\
\label{eq:m19}& \leq &-Cc_0 \int_{\Sigma_\lambda}U_\lambda(y)(\frac{1}{|\tilde x -y|^{n+\alpha}}-\frac{1}{|\tilde x -y^\lambda|^{n+\alpha}})dy.
\end{eqnarray}
Due to the fact that $U_\lambda(x)\not \equiv0$ in $\Sigma_\lambda\backslash \Omega$, \eqref{eq:m19} implies that $${\cal F}_\alpha (u_\lambda(\tilde x))-{\cal F}_\alpha (u(\tilde x))+c_{1 }(\tilde x)U_\lambda(\tilde x)+c_2(\tilde x)V_\lambda(\tilde x)<0.$$
This   contradicts with \eqref{eq:m6}. So $U_\lambda(x)  \equiv0 $ in $\Sigma_\lambda$.
Together with $U_\lambda(x^\lambda)=-U_\lambda(x),$  we arrive at
$$U_\lambda(x)\equiv0,~x\in  {\R}^n.$$
Similarly, one can show that if  $U_\lambda(x)$ attains 0 at one point in $\Omega$, then both $U_\lambda(x)$ and $V_\lambda(x)$ are identically 0 in $ {\R}^n$.
This completes the proof.
\end{proof}

\begin{thm}\label{thmm3}( decay at infinity  )
Assume that $\Omega$ is a subset of $\Sigma_{\lambda}$ and $U_\lambda(x)\in C_{loc}^{1,1}(\Omega) \cap {\cal L}_{\alpha},~ V_\lambda(x)\in C_{loc}^{1,1}(\Omega) \cap {L}_{2s}$,  $U_\lambda(x)$ and  $V_\lambda(x)$ are lower semi-continuous on $\bar \Omega$.
If $U_\lambda(x)$ and  $V_\lambda(x)$ satisfy
\begin{equation}\label{eq:mb11}
\left\{\begin{array}{ll}
{\cal F}_{\alpha}(u_\lambda(x)) -{\cal F}_{\alpha}(u(x))+c_{1 }(x)U_\lambda(x)  +c_{ 2}(x)V_\lambda(x)\geq0  ,&\mbox x \in \Omega,\\[2mm]
(-\lap)^{2s} V_\lambda(x)   +c_{3}(x)U_\lambda(x) \geq0,  &\mbox x \in \Omega,\\[2mm]
U_\lambda(x),V_\lambda(x) \geq0&\mbox x \in \Sigma_\lambda \backslash\Omega,\\[2mm]
U_\lambda(x^{\lambda})=-U_\lambda(x),
V_\lambda(x^{\lambda})=-V_\lambda(x) , &\mbox x \in \Sigma_\lambda,
\end{array}
\right.
\end{equation}
where
\begin{equation}\label{eq:mbb11}
c_{ 1}(x),c_{ 2}(x)\sim o(\frac{1}{|x|^\alpha}),~ c_{3}(x)\sim o(\frac{1}{|x|^{2s}}), ~\mbox{for}~ |x|~\mbox{large} ,\end{equation}
and $$ c_{ 2}(x),~c_{3}(x)< 0.$$
Then there exists a constant $R_0>0$ such that if $$U_\lambda(\tilde x )=\underset {\Omega}{\min}  \, U_\lambda(x)<0,~~V_\lambda(\bar x )=\underset{\Omega}{\min}\, V_\lambda(x)<0 ,$$
then \begin{equation}\label{eq:mbbb11} |\tilde x|\leq R_0~\mbox{or}~|\bar x|\leq R_0.\end{equation}
\end{thm}
\begin{rem} The $x_1$ direction can be chosen arbitrarily, whereas the domain $\Sigma_\lambda$
changes correspondingly, hence the results \eqref{eq:mbbb11} also hold when the problem is set in another direction. \end{rem}

\begin{proof}  By the assumptions,  there exists  $\tilde x \in \Omega,$
such that $$U_\lambda(\tilde x)=\underset{\Omega}{\min}U_\lambda(x)<0.$$
Direct calculation shows that
$$
\aligned {\cal F}_{\alpha}(u_\lambda(\tilde x)) -{\cal F}_{\alpha}(u(\tilde x))
=&C_{n,\alpha}\mbox{ P.V. } \int_{\Sigma_\lambda} \frac{F(u_\lambda(\tilde x)-u_\lambda(y))-F(u(\tilde x)-u(y))}{|\tilde x-y|^{n+\alpha}}\\+&C_{n,\alpha}\mbox{ P.V. } \int_{ \Sigma_\lambda} \frac{F(u_\lambda(\tilde x)-u(y))-F(u(\tilde x)-u_\lambda(y))}{|\tilde x-y^\lambda|^{n+\alpha}}\\
\leq&C_{n,\alpha}\mbox{ P.V. } \int_{\Sigma_\lambda}\frac{F'(\cdot)2U_\lambda(\tilde x)}{|\tilde x-y^\lambda|^{n+\alpha}}dy\\
\leq&2C_{n,\alpha}c_0U_\lambda(\tilde x)\int_{\Sigma_\lambda}\frac{1}{|\tilde x-y^\lambda|^{n+\alpha}}dy.
\endaligned
$$
For each fixed $\lambda$, there exists  $C>0$  such that for $\tilde x \in \Sigma_\lambda$ and $|\tilde x|$ sufficiently large (see \cite{WY}), the following estimate holds
\begin{equation}\label{eq:mb9}
\int_{\Sigma_\lambda}\frac{1}{|\tilde x-y^\lambda|^{n+\alpha}}dy\geq \int_{(B_{3|\tilde x|}(\tilde x)\backslash B_{2|\tilde x|}(\tilde x))\cap \tilde{ \Sigma}_\lambda} \frac{1}{|\tilde x-y|^{n+\alpha}}dy \sim \frac{C}{|\tilde x|^\alpha}.
\end{equation}
Hence
\begin{equation}\label{eq:mb10}
{\cal F}_{\alpha}(u_\lambda(\tilde x)) - {\cal F}_{\alpha}(u(\tilde x))+c_{1 }(\tilde x)U_\lambda(\tilde x) \leq \frac{CU_\lambda(\tilde x)}{|\tilde x|^\alpha}<0.
\end{equation}
Combining \eqref{eq:mb10} with \eqref{eq:mb11}, it is easy to deduce
 \begin{equation}\label{eq:mb12}
 V_\lambda(\tilde x)<0,
 \end{equation}
and
\begin{equation}\label{eq:mb13}
 U_\lambda(\tilde x)\geq -Cc_{ 2}(\tilde x)|\tilde x|^\alpha V_\lambda(\tilde x).
 \end{equation}
Using \eqref{eq:mb12}, there exists  $\bar x$ such that
$$V_\lambda(\bar x)=\underset{\Omega}{\min}V_\lambda(x)<0.$$
Similar to the derivation of \eqref{eq:m20}, we can derive
\begin{equation}\label{eq:mbb13}(-\lap)^sV_\lambda(\bar x)) \leq \frac{CV_\lambda(\bar x)}{|\bar x|^{2s}}<0.\end{equation}
Combing \eqref{eq:mb11}, \eqref{eq:mb13} and \eqref{eq:mbb13}, we have for $\lambda$ sufficiently negative, $$
\aligned
0&\leq  (-\lap)^sV_\lambda  (\bar x)  +c_{3}(\bar x)U_\lambda(\bar x) \\
&\leq \frac{CV(\bar x)}{|\bar x|^{2s}}+c_{3}(\bar x)U_\lambda(\tilde x)\\
&\leq C(\frac{V_\lambda(\bar x)}{|\bar x|^{2s}}-c_{3}(\bar x)c_{ 2}(\tilde x)|\tilde x|^\alpha V_\lambda(\tilde x))\\
&\leq C(\frac{V_\lambda(\bar x)}{|\bar x|^{2s}}-c_{3}(\bar x)c_{ 2}(\tilde x)|\tilde x|^\alpha V_\lambda(\bar x))\\
&\leq\frac{CV(\bar x)}{|\bar x|^{2s}}(1-c_{ 2}(\tilde x)|\tilde x|^\alpha c_{3}(\bar x)|\bar  x|^{2s}),
\endaligned
$$
which shows that $1\leq c_{ 2}(\tilde x)|\tilde x|^\alpha c_{3}(\bar x)|\bar  x|^{2s}.$
However,  from  \eqref{eq:mbb11} we have $$ c_{ 2}(\tilde x)|\tilde x|^\alpha c_{3}(\bar x)|\bar  x|^{2s}<1 $$ for  $|\tilde x|$ and $|\bar x|$ sufficiently large.
This contradiction explains that \eqref{eq:mbbb11} must be true.
This completes the proof.
\end{proof}

\section{Proof of the main result}
This section is contributed in proving Theorem \ref{thmm}, or in other words in obtaining radial symmetry of positive solutions to \eqref{eq:m1}. Since the equivalence of problems \eqref{eq:m1} and \eqref{eq:m3}, we only need to prove Theorem \ref{thmm1}.

{\bf Proof of Theorem \ref{thmm1}.}  Choose an arbitrary direction as the $x_1$-axis, the proof is divided into two steps.

Step 1.  \textit{Start moving the plane $T_\lambda$ from $-\infty$ to the right in $x_1$-direction.}

We will show that for $\lambda$ sufficiently negative,
\begin{equation}\label{eq:mb16}U_\lambda(x)\geq0,~V_\lambda(x)\geq0,~\forall~x\in \Sigma_\lambda.\end{equation}

If  \eqref{eq:mb16} is violated,  then   there are the following 3 possibilities:

$(a)$ both $U_\lambda $  and $V_\lambda  $ are negative in some subsets of $  \Sigma_\lambda$; or

$(b)$ only  $U_\lambda  $ is negative in a subset of $\Sigma_\lambda$; or

$(c)$ only  $V_\lambda  $ is negative in a subset of $\Sigma_\lambda$.

In order to apply Theorem \ref{thmm3}, we first need to rule out possibilities $(b)$ and $(c)$. Then we can prove that the case $(a)$ will not happen. Hence, \eqref{eq:mb16} is true.

 Now we prove that $(b)$ is impossible. If not, we assume that $U_\lambda $ is negative at some point in $\Sigma_\lambda$.
We have $ U_\lambda(x)=V_\lambda(x)\equiv 0, x\in T_\lambda$.
For the fixed $\lambda$, the assumption \eqref{eq:m4} implies that
$$u(x)\rightarrow0,~\mbox{as}~|x|\rightarrow+\infty.$$
Since $|x^\lambda|\rightarrow+\infty$, as $|x|\rightarrow+\infty$,  it follows
$$u_\lambda(x)=u(x^\lambda)\rightarrow0. $$
Thus we have
\begin{equation}\label{eq:maa17}U_\lambda(x)\rightarrow 0,~\mbox{as}~|x|\rightarrow+\infty.\end{equation}
Similarly, one can show that
\begin{equation}\label{eq:aa117}V_\lambda(x)\rightarrow 0,~\mbox{as}~|x|\rightarrow+\infty.\end{equation}
Therefore there exists an $\tilde x \in \Sigma_\lambda$ such that $$U_\lambda(\tilde x)=\underset{\Sigma_\lambda}{\min}U_\lambda(x)<0.$$
The same estimates as in the proof of Theorem \ref{thmm3} yield that
\begin{equation}\label{eq:mbm10}
{\cal F}_{\alpha}(u_\lambda(\tilde x)) - {\cal F}_{\alpha}(u(\tilde x))+c_{1 }(\tilde x)U_\lambda(\tilde x) \leq \frac{CU_\lambda(\tilde x)}{|\tilde x|^\alpha}<0 .
\end{equation}
However combining \eqref{eq:m3} with the mean value theorem, we obtain
\begin{equation}
\label{eq:mbmm019}\aligned {\cal F}_\alpha (u_\lambda(\tilde x))-{\cal F}_\alpha (u(\tilde x))
+(\omega-v(\tilde x)r\xi^{r-1})U_\lambda(\tilde x)= u_\lambda^r(\tilde x)V_\lambda(\tilde x) ,\endaligned
\end{equation}
where $\xi$ is a value between $u_\lambda(\tilde x)$ and $u(\tilde x)$.
From  \eqref{eq:mbmm019} and \eqref{eq:mbm10} with $c_1(\tilde x)=(\omega-v(\tilde x)r\xi^{r-1})$, we get $ V_\lambda(\tilde x)<0$, which contradicts to $V_\lambda(x)\geq 0$. Hence case
$(b)$ cannot happen. By similar discussion,  one can rule out case $(c)$.

Next we will prove that the case $(a)$ will not happen. In this case, $U_\lambda$ and $V_\lambda$ have negative minimum points separately.
So we can conclude there exists an $\bar x \in \Sigma_\lambda$ such that \begin{equation}\label{eq:mww3003}V_\lambda(\bar x)=\underset{\Sigma_\lambda}{\min}V_\lambda(x)<0.\end{equation}
Furthermore, we claim that \begin{equation}\label{eq:mmw10}
 U_\lambda(\bar x)<0.
\end{equation}
If not, it follows from \eqref{eq:m3} and the mean value theorem that $$ (-\lap)^sV_\lambda(\bar x) =q\eta^{q-1}U_\lambda(\bar x)\geq 0,$$
where   $\eta$ is valued between $u_\lambda(\bar x)$ and $u(\bar x)$, which contradicts with \eqref{eq:mbb13}.
  This contradiction deduces \eqref{eq:mmw10}.

Therefore, via the above conclusions and the mean value theorem, we arrive at
\begin{equation}\label{eq:mmm019}\left\{\begin{array}{ll}
{\cal F}_\alpha (u_\lambda(\tilde x))-{\cal F}_\alpha (u(\tilde x))+(\omega-v(\tilde x)ru^{r-1} (\tilde x))U_\lambda(\tilde x)-u ^r(\tilde x)V_\lambda(\tilde x)\geq 0,\\[2mm]
(-\lap)^sV_\lambda(\bar x)  -qu^{q-1}  (\bar x)U_\lambda(\bar x)\geq 0.
 \end{array}
\right.\end{equation}
By Theorem \ref{thmm3}, it suffices to check the decay rate at the points   where $V_\lambda(x)$ and $U_\lambda(x)$ are negative respectively.
At those points  for $|x|$ sufficiently large, the decay assumptions  \eqref{eq:m4} and  \eqref{eq:m5} instantly yields that
\begin{equation}\label{eq:mw33} \aligned &c_1(\tilde x)= \omega-v(\tilde x)ru^{r-1} (\tilde x)\sim o(\frac{1}{|\tilde x|^\alpha}),~c_2(\tilde x)=-u ^r(\tilde x)\sim o(\frac{1}{|\tilde x|^{\alpha}}), \\   & c_3(\bar x)=-qu^{q-1}  (\bar x)\sim o(\frac{1}{|\bar x|^{2s}}).\endaligned \end{equation}
Consequently, there exists $R_0>0,$    and  it holds by Theorem \ref{thmm3} that
\begin{equation}\label{eq:maa16}|\tilde x|\leq R_0 ~\mbox{or}~|\bar x|\leq R_0.\end{equation}
Without loss of generality, we may assume
\begin{equation}\label{eq:maa18}|\tilde x|\leq R_0 .\end{equation}
 Hence we have for $\lambda$ sufficiently negative, \begin{equation}\label{eq:bb108}U_\lambda(x)\geq0,~\forall ~x\in\Sigma_\lambda.\end{equation}
Now we claim  $V_\lambda(x)\geq 0$ in $\Sigma_\lambda$. Otherwise, we obtain \eqref{eq:mww3003}
 and then   it admits from \eqref{eq:mbb13} that
 \begin{equation}\label{eq:mbbb13}(-\lap)^sV_\lambda(\bar x)  \leq \frac{CV_\lambda(\bar x)}{|\bar x|^{2s}}< 0.\end{equation}
However, the second equation of  \eqref{eq:mmm019} with   $u(x),v(x)>0$ yields $$(-\lap)^sV_\lambda(\bar x) \geq qu^{q-1}  (\bar x)U_\lambda(\bar x)\geq 0.$$  This contradicts to  \eqref{eq:mbbb13}, which means that  $V_\lambda(x)$ is nonnegative  $\Sigma_\lambda$. So (a) will not happen. Therefore \eqref{eq:mb16} is proved.

Step 2. \textit{Keep moving  the planes to the right till the limiting position $T_{\lambda_0}$ as long as
\eqref{eq:mb16} holds}.

Let $$\lambda_0=\sup\{\lambda \mid  U_\mu(x),~ V_\mu(x) \geq 0, ~x\in \Sigma_\mu, ~\mu\leq \lambda\},$$ then the behaviors of  $ u$ and $v$ at infinity guarantee $\lambda_0<\infty$.

In this part, we show that
 \begin{equation}\label{eq:mb61}U_{\lambda_0}(x)\equiv0, ~ V_{\lambda_0}(x) \equiv 0,~ x\in \Sigma_{\lambda_0}. \end{equation}
By the definition of $\lambda_0$ and $(iii)$ of Theorem \ref{eq:m2}, we first point that either \eqref{eq:mb61} or
\begin{equation}\label{eq:mm66}U_{\lambda_0}(x)>0, ~ V_{\lambda_0}(x) >0,~ \forall~ x\in \Sigma_{\lambda_0} \end{equation}
holds.

In fact, if  \eqref{eq:mb61} is violated, then \eqref{eq:mm66} must be true. In this case, it can be shown that one can move the plane $T_\lambda$  further to the right such that \eqref{eq:mb16}  is still valid. More precisely, our remaining task is to prove that  there exists a small $\epsilon>0,$ such that for any
$\lambda\in (\lambda_0,\lambda_0+\epsilon)$, it holds
\begin{equation}
\label{eq:mb17}U_\lambda(x)\geq 0,~V_\lambda(x)\geq 0, ~\forall~ x\in \Sigma_\lambda,
\end{equation}
which is a contradiction to the definition of $\lambda_0$. Hence \eqref{eq:mb61} must be true.

The proof of \eqref{eq:mb17} by using Theorem \ref{thmm2} and Theorem \ref{thmm3} is given in the following.
Suppose that \eqref{eq:mb17} is false,    then due to the argument
in Step 1.  both $U_\lambda$ and  $V_\lambda$ achieve their negative minima in $\Sigma_\lambda$, $i.e.$ ($b$) and ($c$) are impossible. The proof of this conclusion is the same as step 1, here we omit. Then we will derive
contradictions by showing that these minima can fall nowhere in $\Sigma_\lambda$, that is, we can rule out case ($a$) and so \eqref{eq:mb17} is proved.

 Next to prove  \eqref{eq:mb17}.  In case ($a$), let $\tilde x$ and $\bar x$ be the minimum points of $U_\lambda$ and $V_\lambda$ separately, i.e.
$$
U_\lambda(\tilde x )=\underset {\Sigma_\lambda}{\min} U_\lambda(x)<0,~~V_\lambda(\bar x )=\underset{\Sigma_\lambda}{\min} V_\lambda(x)<0.
$$

Let $R_0$ be determined in Theorem \ref{thmm3}. From \eqref{eq:mm66}, for any  $\delta>\epsilon>0$, we have
$$U_{\lambda_0}(x)\geq c_0>0,~V_{\lambda_0}(x)\geq c_0>0,~  \forall ~ x\in \overline{\Sigma_{\lambda_0-\delta}\cap B_{R_0}(0)} .$$
Using the continuity of $U_\lambda(x)$ and $ V_\lambda(x)$ with respect to $\lambda$, there exists   $\epsilon>0,$  such that  for all  $\lambda\in (\lambda_0,\lambda_0+\epsilon),$ it holds
\begin{equation}\label{eq:b19}
U_\lambda(x)\geq 0,~V_\lambda(x)\geq 0, ~\forall~ x\in \overline{\Sigma_{\lambda_0-\delta}\cap B_{R_0}(0)} .
\end{equation}
Now we are ready to consider the following three cases:

\textbf{Case 1}. $\tilde x \in B_{R_0}(0)\bigcap (\Sigma_{\lambda_0+\epsilon}\backslash \Sigma_{\lambda_0-\delta})$ and $\bar x \in \Sigma_\lambda \bigcap B^c_{R_0}(0).$

 Similar to the derivation of \eqref{eq:mb13}, we have
 \begin{equation}\label{eq:bb8}
 U_\lambda(\tilde x)\geq -Cc_{ 2}(\tilde x)l^\alpha V_\lambda(\tilde x)
\end{equation}
and
$$\aligned
0&\leq (-\lap)^s V_\lambda(\bar x)  +c_{3}(\bar x)U_\lambda(\bar x) \\
&\leq\frac{CV_\lambda(\bar x)}{|\bar x|^{2s}}+c_{3}(\bar x)U_\lambda(\tilde x)\\
&\leq C\{\frac{V_\lambda(\bar x)}{|\bar x|^{2s}}-c_{3}(\bar x)c_{ 2}(\tilde x)l^\alpha V_\lambda(\tilde x)\}\\
&\leq C\{\frac{V_\lambda(\bar x)}{|\bar x|^{2s}}-c_{3}(\bar x)c_{ 2}(\tilde x)l^\alpha V_\lambda(\bar x)\}\\
&\leq C \frac{V_\lambda(\bar x)}{|\bar x|^{2s}}[1-c_{ 2}(\tilde x)l^{\alpha}c_{3}(\bar x)|\bar x|^{2s}].
\endaligned$$
Hence \begin{equation}\label{eq:mbbb20}1\leq c_{ 2}(\tilde x)l^{\alpha}c_{3}(\bar x)|\bar x|^{2s}.\end{equation}
However, we know by  \eqref{eq:mw33} that $|c_{3}(\bar x)|\bar x|^{2s}|$ is small for $|\bar x|$ sufficiently large. Due to the facts that $l=\epsilon+\delta$ is small and $c_{ 2}(\tilde x)$ is bounded from below in
$\Sigma_{\lambda_0+\epsilon}\backslash \Sigma_{\lambda_0-\delta}$, we obtain that $|c_{ 2}(\tilde x)l^{\alpha}|$ is small. Consequently,  we have that $c_{ 2}(\tilde x)l^{\alpha}c_{3}(\bar x)|\bar x|^{2s}<1$, which contradicts with \eqref{eq:mbbb20}. Therefore in case 1, (a) will not happen, hereby  \eqref{eq:mb17} is proved.

\textbf{Case 2}. $\bar x \in B_{R_0}(0)\bigcap (\Sigma_{\lambda_0+\epsilon}\backslash \Sigma_{\lambda_0-\delta})$ and   $\tilde x \in \Sigma_\lambda \bigcap B^c_{R_0}(0).$

The validity of \eqref{eq:mb17} can be proved similarly to the discussion in \textbf{Case   1}, which is omitted here.

\textbf{Case 3}.  $\tilde x, \bar x \in B_{R_0}(0) \bigcap (\Sigma_{\lambda_0+\epsilon}\backslash \Sigma_{\lambda_0-\delta})$.

By taking $l=\delta+\epsilon$ in \eqref{eq:mb7}, together with \eqref{eq:mmm019}, we arrive at
\begin{equation}\label{eq:mbb8}
U_\lambda(\tilde x)\geq -Cc_{ 2}(\tilde x)l^\alpha V_\lambda(\tilde x).
\end{equation}
From \eqref{eq:m99}, we  derive
$$(-\lap)^sV_\lambda(\bar x)  \leq \frac{CV_\lambda(\bar x)}{l^{2s}}  <0.$$
Noting \eqref{eq:mbb8}, we have for $l$ sufficiently small,
$$\aligned
0&\leq  (-\lap)^sV_\lambda(\bar x)  +c_{3}(\bar x)U_\lambda(\bar x) \\
&\leq\frac{CV_\lambda(\bar x)}{l^{2s}}+c_{3}(\bar x)U_\lambda(\tilde x)\\
&\leq C\{\frac{V_\lambda(\bar x)}{l^{2s}}-c_{3}(\bar x)c_{ 2}(\tilde x)l^\alpha V_\lambda(\tilde x)\}\\
&\leq C\{\frac{V_\lambda(\bar x)}{l^{2s}}-c_{3}(\bar x)c_{ 2}(\tilde x)l^\alpha V_\lambda(\bar x)\}\\
&\leq C \frac{V_\lambda(\bar x)}{l^{2s}}[1-c_{ 2}(\tilde x)c_{3}(\bar x)l^{\alpha+{2s}}]<0.
\endaligned$$
This contradiction shows that $(a)$ does not happen. Therefore the statement \eqref{eq:mb17} is correct.

Now we have shown that $ U_{\lambda_0}(x)\equiv 0,~V_{\lambda_0}(x)\equiv0, ~ x\in \Sigma_{\lambda_0}.$
Since the $x_1$-direction can be chosen arbitrarily, we have proven that $u(x)$  must be radially symmetric  about some point  in $\R^n$.
Also the monotonicity follows easily from the argument.
This completes the proof of  Theorem \ref{thmm1}.

\section*{Acknowledgments}
The work was carried out when the first author visits University of Mannheim in Germany. Pengyan Wang is supported by the scholarship of NPU's  exchange funding program. Li Chen is partially supported by DFG Project CH 955/4-1. Pengcheng Niu is supported by National Natural Science Foundation of China (Grant No.11771354).

\end{document}